\documentstyle[12pt]{article}

\newcommand{\sect}[1]{\setcounter{equation}{0}\section{#1}}
\newcommand{\subsect}[1]{\subsection{#1}}

\def\be{\begin{equation}}
\def\ee{\end{equation}}
\def\bea{\begin{eqnarray}}
\def\eea{\end{eqnarray}}

\def\aap{a_+}
\def\aam{a_-}
\def\baap{\overline{ a}_+}
\def\baam{\overline{ a}_-}

\def\bbp{b_+}
\def\bbm{b_-}

\def\bp{\beta_+}
\def\bm{\beta_-}

\def\aa{N}
\def\baa{\overline{N}} 
 
\def\bb{M}

\def\ap{A_+}
\def\am{A_-}
\def\bap{\overline{A}_+}
\def\bam{\overline{A}_-}

\def\jj{J_3} 
\def\jp{J_+} 
\def\jm{J_-}
\def\bjj{\overline{J}_3} 
\def\bjp{\overline{J}_+} 
\def\bjm{\overline{J}_-}

\def\>#1{{\bf #1}}                 

\def\Z{\bf Z} 
\def\1{\'{\i}}                           
\def\R{{\rm I\kern-.2em R}}

\begin{document}

\thispagestyle{empty}

 \hfill q-alg/9611031
\bigskip\bigskip
\bigskip

\begin{center} 

{\LARGE{\bf{Boson representations,}}} 

{\LARGE{\bf{non-standard quantum algebras}}} 

{\LARGE{\bf{and contractions}}} 

\end{center}

\bigskip\bigskip

\begin{center} Angel Ballesteros$^\dagger$, Francisco J.
Herranz$^\dagger$ and Javier Negro$^\ddagger$
\end{center}

\begin{center} {\it {  $^\dagger$ Departamento de F\1sica, Universidad
de Burgos} \\   Pza. Misael Ba\~nuelos, 
E-09001, Burgos, Spain}
\end{center}

\begin{center} {\it {  $^\ddagger$ Departamento  de F\1sica Te\'orica,
Universidad de Valladolid } \\   E-47011, Valladolid, Spain}
\end{center}

\bigskip\bigskip\bigskip

\begin{abstract} 
A Gelfan'd--Dyson mapping is used to generate   a one-boson realization for
the non-standard quantum deformation of  $sl(2,\R)$ which directly provides
its infinite and finite dimensional irreducible representations. Tensor
product decompositions are worked out for some examples.  Relations between
contraction methods and boson realizations are also explored in several
contexts. So, a class of two-boson representations for the non-standard
deformation of 
$sl(2,\R)$ is introduced and contracted to the non-standard quantum (1+1)
Poincar\'e representations. Likewise, a  quantum extended Hopf $sl(2,\R)$
algebra is constructed and the Jordanian $q$-oscillator algebra
representations are obtained from it by  means of another contraction
procedure.
\end{abstract}

\vspace{3cm}

\newpage

%%%%%%%%%%%%%%%%% introduccion %%%%%%%%%%%

\sect {Introduction}

Boson realizations of many symmetry algebras and superalgebras are known
to be useful in many problems of  Condensed Matter \cite{CM} and Nuclear
Physics \cite{NP}. Among the existent bosonization processes, we shall fix our
attention in the so-called Gelfan'd--Dyson (GD) mapping of  $sl(2,\R)$
\cite{dy}, initially introduced in spin systems. The aim of this paper is to
show that deformed GD type realizations are the most appropriate
tools in order to construct the representation theory of non-standard quantum
$sl(2,\R)$ and other non-standard quantum algebras linked to it by means of
contraction limits. Therefore, we hope that the results reached in this paper
can be directly applied to deformed shell models or coherent states methods
where GD maps have been proven very successful.

Firstly, we recall  that the standard deformation of $sl(2,\R)$
\cite{Ku,Dr,Ji} is associated to the (constant) solution
$r=J_+\wedge J_-$ of the modified classical Yang--Baxter
equation (YBE). This quantum algebra has been fully developed
and extensively applied (see \cite{PC} ). However, there also exists a
non-standard deformation linked to the solution $r=J_3\wedge J_+$ of the  
classical YBE. This deformation (which was firstly introduced at a quantum
group level \cite{Manin,Zakr}, and later as a quantum Hopf algebra \cite{Ohn})
has recently attracted  much attention. For instance, it has been applied to
build up higher dimensional non-standard quantum algebras \cite{Beyond} as
well as the non-standard
$q$-differential calculus \cite{Karim,azca}. Its universal quantum $R$-matrix 
\cite{Iran,rr} and its irreducible representations \cite{dob,ab,abc}  have been
also studied.

Furthermore, it is rather remarkable that  there exists a close
relationship between $U_zsl(2,\R)$ and the non-standard quantum (1+1)
Poincar\'e
\cite{Iranc,tt} and oscillator algebras  \cite{osc}. In particular, all of
them have a similar Hopf subalgebra determined by the two generators involved
in the classical $r$-matrix. As an important  consequence, they have a
formally identical universal $R$-matrix. We will show that most of these
common features can be explained by a contraction  scheme  connecting all
these non-standard quantum algebras.

As we shall see in section 2, the full
representation theory of $sl(2,\R)$ can be straightforwardly derived with the
aid of the one-boson GD realization. Afterwards, we build the one-boson
infinite dimensional representations for $U_zsl(2,\R)$ by following the same
appoach. It turns out that their explicit form is somewhat more complicated
than those of the standard deformation \cite{bi,mac} in the sense that they
cannot be obtained by the mere substitutions of numbers by $q$-numbers. The
corresponding finite dimensional representations are deduced in a very
natural way obtaining closed expressions for their matrix elements in any
dimension. It is remarkable that these results are in plain agreement with
those derived in \cite{ab,abc} by using  a recurrence method in another basis.
The main advantage of the basis here used is the very simple form that the
quantum universal $R$-matrix presents; this fact,  combined with the
representation theory just derived, is used to provide some explicit
$R$-matrices.
 
Although these deformed representations may appear  ackward to handle, we
know that they should be fully consistent with the deformed composition of
representations given by the quantum coproduct. In section 3 we illustrate 
such a problem  for  some low dimensional representations  showing that
complete reducibility holds and preserves the same well
known classical angular momentum decomposition rules in tensor product spaces.
However, non-standard Clebsch-Gordan coefficients are shown to be essentially
different to those of the standard deformation.

The remaining sections of the paper are devoted to the study of the relations
between boson realizations and contractions in classical and deformed
frameworks. In particular, two-boson GD representations for
classical and quantum $sl(2,\R)$ are introduced  in section 4. They are shown
to be the most adequate objects to obtain the representations of the
non-standard quantum Poincar\'e algebra by means of a contraction process. We
also present in section 5 a suitable quantum deformation of the
(pseudo)extended Lie algebra $sl(2,\R)$ together with its one-boson
representations so that they give rise, also through a contraction procedure,
to the representations of the non-standard quantum oscillator algebra (also
quoted as the Jordanian
$q$-oscillator). The
 extension here introduced contains some
interesting features that will be discussed. In  
  section 6  the non-standard  $sl(2,\R)$,
Poincar\'e and oscillator algebras are presented as quadratic structures  by
using a deformed boson algebra. Finally, some remarks end the paper.

%%%%%%%%%%%%%%%%%%%%%%%%%%%%%%%%%%%%%%%%%%%%%%%%%
%%%%%%%%%%%%%%%%% One-boson %%%%%%%%%%%%%%%%%%%%%%%%

\sect {One-boson $U_zsl(2,\R)$ representations} 

\subsect{Classical one-boson representations}

To start with we shall deal with the classical Lie algebra $sl(2,\R)$ whose
generators
$\{\jj,\jp,\jm\}$ obey the commutation rules
\be
[\jj,\jp]=2\jp,\qquad  [\jj,\jm]=-2\jm,\qquad [\jp,\jm]=\jj .
\label{ca}
\ee
This algebra is isomorphic to $so(2,1)$, the Lie algebra generating the group
of motions of a (1+1) De Sitter space with non-zero constant curvature, where 
$\jj$  generates the  boosts and $\jp$, $\jm$ translations along the
light-cone. An alternative physical interpretation for  $sl(2,\R)$ is to
consider it as the infinitesimal generators of the conformal group of a 
one-dimensional space; in this sense, $\jj$ would generate dilations, $\jp$
translations and $\jm$ special conformal transformations. Obviously, these
different interpretations come from different  representations of the algebra
on the spaces linked to the physical problem under consideration. 
 
The irreducible representations of $sl(2,\R)$ are characterized
by the eigenvalue of the quadratic Casimir element 
\be
{\cal C}=\frac 12 \jj^2 +\jp\jm+\jm\jp .
\label{cb}
\ee
If the generators $\{\aam,\aap\}$ close a boson algebra, i.e.,
$[\aam,\aap] = 1$, then, the realization of   $sl(2,\R)$ given by
\be 
 \jp=\aap,\qquad
 \jj=2\aap\aam +\beta 1,\qquad
 \jm=-\aap\aam^2-\beta\aam , 
\label{cd}
\ee 
where $\beta$ is a free parameter, is known as the GD
one-boson realization \cite{dy,agi}. Now  we will see how the GD map
(\ref{cd}) can be used in order to get easily any of the $sl(2,\R) \approx
su(1,1)$ irreducible representation series \cite{vk}.

\smallskip

\noindent
{\it (i) Lower bounded representations.} When the operators
$\aam$, $\aap$ act in the usual way on the number states Hilbert space
spanned by $\{|m\rangle\}_{m=0}^{\infty}$, i.e.,
\be
\aap |m\rangle = \sqrt{m+1}\,|m+1\rangle, \qquad 
\aam |m\rangle = \sqrt{m}\,|m-1\rangle ,\label{cdd}
\ee
(\ref{cd}) leads to a lower bounded representation:
\bea 
 &&\jp|m\rangle =\sqrt{m+1}\, |m+1\rangle ,\qquad
  \jj|m\rangle =(2m+\beta)\, |m\rangle,\cr
 && \jm|m\rangle =-\sqrt{m}\,(m-1+\beta)\, |m-1\rangle .
\label{cf}
\eea 
The casimir eigenvalue being
\be
{\cal C}=\beta(\beta/2 -1).
\label{cg}
\ee
For negative integer values of $\beta$, hereafter  denoted as  $\bm\in{\bf
Z}^-$, the representation (\ref{cf}) is reducible leading to a finite
dimensional irreducible quotient representation of dimension $|\bm -1|$. For
instance, $\bm=-1$ (${\cal C}=3/2$) provides the two-dimensional
representations of
$sl(2,\R)$ by setting 
$|2\rangle\equiv 0$:
\be 
 J_+|0\rangle=|1\rangle\quad
J_+|1\rangle=0\quad
J_3|0\rangle=- |0\rangle\quad
J_3|1\rangle= |1\rangle\quad
 J_-|0\rangle=0\quad
J_-|1\rangle= |0\rangle .
\label{cgb}
\ee 
The numbers $\langle m|X|m'\rangle$ where $\langle m|m'\rangle=\delta_{m,m'}$
give the matrix elements of these representations;  in  the previous example,
we have 
\be
\jp=\left(\begin{array}{ll}
. & .   \cr 
1  & .  \cr  
\end{array}\right) ,
\qquad
\jm=\left(\begin{array}{ll}
. & 1   \cr 
.  & .  \cr  
\end{array}\right) ,\qquad
\jj=\left(\begin{array}{ll}
-1 & .   \cr 
. & 1  \cr  
\end{array}\right)  .
\label{cgc}
\ee

With this notation in mind, we shall schematically display the matrix form
that the aforementioned representations take as a future
reference with respect to the deformed algebra:
\bea
&&\jp=\left(\begin{array}{llllll}
0 & . & . & . & . & \cr 
1  & 0 & . & . & . &\cr . 
& {\sqrt{2}} & 0 & . & . & \cr
. & . & {\sqrt{3}} & 0 & . &\cr 
. & . & . &{\sqrt{4}} & 0 &\cr
  &    &   &  &   &\ddots\cr
\end{array}\right),\label{ci}\\
&&\jj=\left(\begin{array}{llllll}
\beta & .   & . & . & . &\cr . &  2+\beta &. & . & . &
   \cr . & . & 4+\beta & . & .& \cr 
. & . & . &6+\beta & . & \cr   
.&.   &  . & .  & 8+\beta  & \cr
&   &   &   &   &\ddots
\end{array}\right) ,\label{cj}\\
&&\jm=-\left(\begin{array}{llllll}
0 & \beta   & . & . & . &\cr . & 0 & {\sqrt{2}}(1+\beta) & . & . &
   \cr . & . & 0 & {\sqrt{3}}(2+\beta) & .& \cr 
. & . & . &0 & {\sqrt{4}}(3+\beta) &
   \cr 
. & . & . &. & 0 &
   \cr 
  &   &   &   &   &\ddots
\end{array}\right) .
\label{ck}
\eea

\smallskip

\noindent
{\it (ii) Upper bounded representations.}  Quite similar upper bounded
representations can be defined in the suplementary space
$\{|m\rangle\}_{m=-\infty}^{-1}$. However, in order to avoid the complex
numbers in the accompanying square roots (\ref{cf}) inside this space, we
shall redefine the basis vectors in the form
$|m\rangle \to \frac{-1}{\sqrt m}|m\rangle$, so that  the boson operators act
as
\be
\aap |m\rangle = -(m+1)\,|m+1\rangle, \qquad 
\aam |m\rangle = -|m-1\rangle ,\label{cka}
\ee
leading to the $sl(2,\R)$ action
\be 
\jp|m\rangle =-(m+1)\, |m+1\rangle ,\quad
\jj|m\rangle =(2m+\bp)\, |m\rangle,\quad
\jm|m\rangle =(m-1+\bp)\, |m-1\rangle .
\label{ckb}
\ee
The finite dimensional representations are now originated for 
$\bp-2\in{\bf Z}^+$, with dimension
$\bp-1$. However, note that in this case the action (\ref{ckb}) allows for
an invariant subspace, so that it is not necessary  to make use of the
quotient mechanism to reach irreducibility. 

These representations are particularly well
suited for describing the differential version of the GD map (\ref{cd}),
\be 
 \jp=\partial_x,\qquad
 \jj=-2x\partial_x +\beta-2 ,\qquad
 \jm=-x^2\partial_x + (\beta-2)x .
\label{ce}
\ee 
The basis functions will be the positive integer powers 
$\{x^n\}_{n=0}^{+\infty}$, with the identification $x^n \equiv
|{-n-1}\rangle$. In particular for the values of the label $\beta$
given by $\bp-2\in {\bf Z}^+$ the support space for the finite
$(\bp - 1)$-dimensional representations is generated by the monomials
$\{ 1,x,x^2,\ldots , x^{\bp-2}\}$.  Note also that (\ref{ce})
reproduces for $\beta=2$ the usual differential realization  of the Lie
algebra of the conformal group for the one-dimensional Euclidean space.

The finite dimensional representations obtained
from the lower (labeled by $\bm$) or upper bounded ones (denoted by
$\bp$) are equivalent whenever $\bp-2 = |\bm|=-\bm$. Indeed in this case
the Casimir (\ref{cg}) is the same ${\cal C}_{\bm} = 
{\cal C}_{\bp}$. Hereafter we shall introduce the notation  
$|\bm -1| = \bp -1 = 2j+1$, where
$j$ is a half positive integer that should be identified  with the label of
the integer $(2j+1)$-dimensional representations of $sl(2,\R) \approx su(1,1)$
\cite{vk}. Using this notation, the Casimir (\ref{cg}) will  turn into the
more familiar expression ${\cal C} = 2j(j+1)$.

\smallskip

\noindent
{\it (iii) Non-bounded representations.} When the operators  $\aam$, $\aap$
act on the space spanned by
$\{|m\rangle\}_{m=-\infty}^{\infty}$ in the form
\be
\aap |m\rangle = \sqrt{m +\lambda+1}\, |m+1\rangle, \qquad 
\aam |m\rangle = \sqrt{m+\lambda}\, |m-1\rangle ,\label{cgg}
\ee
where $0<Re(\lambda)<1$ we get unbounded representations. These are less
familiar because they cannot implement the hermiticity conditions
$(\aam)^{\dag} = \aap$ inside a Hilbert space. However they can  be
substituted in (\ref{cf}), so that if
$\beta-\lambda
\notin \Z$ we get a family of
$sl(2,\R)$ representations with the same eigenvalue (\ref{cg}) which are not
bounded either. Nevertheless we will not need to give more details
because this paper shall be mainly concerned with the bounded cases.

The principal and complementary series of $sl(2,\R) \approx su(1,1)$
representations can be derived from case {\it (iii)}, while the discrete
series come from the other {\it (i)} and {\it (ii)} cases. The implementation
of the hermiticity conditions can be achieved by a wiseful definition of the
inner product. This is a familiar problem of GD mappings which often is
called `the unitarization process' \cite{NP}. Anyway, we shall  not address
these questions here.

\subsect{Quantum one-boson representations}

The Hopf algebra $U_z sl(2,\R)$ deforming the bialgebra generated by the
classical $r$-matrix $r=zJ_3\wedge J_+$ is characterized by the following
coproduct, counit, antipode and commutation rules
(see \cite{rr}):
\bea
&&  \Delta (\jp) =1 \otimes \jp  + \jp \otimes 1,\cr
&& \Delta (\jj) =1 \otimes \jj  + \jj\otimes
e^{2z \jp }, \label{ea} \\
&& \Delta (\jm) = 1 \otimes \jm + \jm\otimes
e^{2 z \jp }, 
\nonumber
\eea
\be
\epsilon(X) =0,\qquad 
 \mbox{for $X\in \{\jj,\jp,\jm\}$},
 \ee
\be
\gamma(\jp)=-\jp,\qquad \gamma(\jj)=-\jj e^{-2z\jp},\qquad 
\gamma(\jm)=-\jm e^{-2z\jp},
\ee
\be 
[\jj,\jp ]= \frac{e^{2 z \jp} -1  } z,\quad 
[\jj,\jm]=-2 \jm +z \jj^2,\quad  [\jp  ,\jm ]= \jj.
\label{eb} 
\ee 
The quantum Casimir is
\be
{\cal  C}_z=\frac 12 \jj\, e^{-2z\jp}\jj +\frac {1-e^{-2z\jp}}{2z}\,\jm
+\jm\,\frac {1-e^{-2z\jp}}{2z} + e^{-2z\jp}-1,
\label{ec}
\ee
and the  universal $R$-matrix takes the form
\be
{\cal  R}=\exp\{-z \jp\otimes \jj\}\exp\{z \jj\otimes \jp\} .
\label{ed}
\ee

A realization of $U_zsl(2,\R)$ in terms of the boson
algebra $[\aam,\aap]=1$  reads
\bea
&& \jp=\aap, \qquad \jj=\frac{e^{2 z \aap} -1  }{z}\,\aam 
+\beta \, \frac{e^{2 z \aap} +1  }{2},\label{ee}\\
&& \jm=-\frac{e^{2 z \aap} -1  }{2\,z}\,\aam^2 - 
\beta \, \frac{e^{2 z \aap} +1  }{2}\,\aam 
- z \beta^2\,\frac{e^{2 z \aap} -1  }{8}.\nonumber
\eea
The Casimir for this realization takes the same expression (\ref{cg}) than
in the classical case.
The limit $z \to 0$ of (\ref{ee}) gives rise to the
GD realization (\ref{cd}) for $sl(2,\R)$ while the Casimir keeps the
same eigenvalue (\ref{cg}) along the whole process.

The essenial point is that GD like quantum formulas (\ref{ee})  allow us to
compute easily  closed expressions for any representation of this quantum
algebra, (in this respect, compare to the derivation of such representations 
given in \cite{ab,abc}). 

Lower bounded representations can be obtained from (\ref{cdd}),  by taking
into account that  \be e^{2 z \aap}|m\rangle =|m\rangle+   
\sum_{k=1}^\infty \frac{(2z)^k }{k!} \, \sqrt{\frac {(m+k)!}{m!} }
\, |m+k\rangle ,
\label{eegg}
\ee
so that the action of (\ref{ee}) on the states
$\{|m\rangle\}_{m=0}^{\infty}$ for any
$\beta$ is obtained:
\bea
&&\jp|m\rangle =\sqrt{m+1}\, |m+1\rangle ,\cr
&&\jj|m\rangle =(2m+\beta)\, |m\rangle +  \sum_{k=1}^\infty \frac{(2z)^k
}{k!} \, \sqrt{\frac {(m+k)!}{m!} }\left(
\frac{2m}{k+1}+\frac{\beta}{2}\right)\, |m+k\rangle ,\label{eg}\\
&&\jm|m\rangle =-\sqrt{m}\,(m-1+\beta) \, |m-1\rangle\cr
&&\ - \sum_{k=1}^\infty \frac{(2z)^k}{k!}
\sqrt{\frac {(m+k)!}{m!} }\left\{
\frac{m}{\sqrt{m+k}}\left(\frac{m-1}{k+1}+\frac{\beta}{2}\right)
\, |m-1+k\rangle +\frac{z\beta^2}{8}\, |m+k\rangle\right\} .
\nonumber
\eea
The general matrix form for $\jp$ is the classical one (\ref{ci}) and for the
two remaining generators we have
\bea
&&\jj=\left (\begin{array}{llllll}
\beta  & . & . & . & .& \cr 
\beta  z &    2 + \beta  & . & . & .& \cr 
{\sqrt{2}} \beta  {z^2} &
  {\sqrt{2}}  ( 2 + \beta   )  z & 4 +
\beta   & . & . &\cr 
\frac {4}{\sqrt{6}} \beta  {z^3} & 
   2 {\sqrt{6}}  ( {2\over 3} + 
  \frac 12   \beta  )  {z^2} & 
  {\sqrt{3}}  ( 4 + \beta   )  z & 6 + \beta 
    & . &\cr 
\frac {4}{\sqrt{6}} \beta  {z^4} & 
  4 {\sqrt{{2\over 3}}}  ( 1 + \beta   )  
   {z^3} & 2 {\sqrt{3}} 
    ( {8\over 3} + \beta   )  {z^2} & 
   2  ( 6 + \beta   )  z & 8 + \beta & \cr
  &    &   &  &   &\ddots\cr
\end{array} \right),\cr
 &&\jm= \nonumber
\eea
{\footnotesize{\be
- \left (\begin{array}{llllll}
0 \!&  \beta  \!\!& . \!\!& . \!\!& . \!\!&\cr 
 \frac 14  {\beta }^2   {z^2} \!\!& \beta  z    \!\!& 
  {\sqrt{2}}  ( 1 + \beta   )     \!\!& . \!\!& . \!\!&\cr 
\frac 1{ 2 {\sqrt{2}} }{{\beta }^2} {z^3}\!\!&
  ( {\sqrt{2}} \beta +\frac 1{2  {\sqrt{2}}} \beta^2
 )   {z^2}   \!\!&  2  ( 1 + \beta   )  z \!\!& 
{\sqrt{3}}  ( 2 + \beta   )  \!\!& . \!\!& \cr 
\frac 1{\sqrt{6} } \beta^2 {z^4} \!\!& 
   ( 2 {\sqrt{{2\over 3}}} \beta + 
\frac {\sqrt{3}}{2\sqrt{2}}\beta^2 ) z^3 \!\!& 
  (\frac{4}{\sqrt{3}}+ 2\sqrt{3} \beta+\frac{\sqrt{3}}{4} \beta^2 )
  z^2 \!\!&  3  ( 2 + \beta   )  z \!\!& \sqrt{4}  ( 3 + \beta   )  \!\!&\cr 
 \frac 1{\sqrt{6} } \beta^2 {z^5}  \!\!& 
 (  2 {\sqrt{{2\over 3}}} \beta +
   {\sqrt{{2\over 3}}} {{\beta }^2} ) {z^4} \!\!& 
   ( \frac{4}{\sqrt{3}} +  \frac{8}{\sqrt{3}} \beta+
\frac{\sqrt{3}}{2} {{\beta }^2} ) {z^3} 
 \!\!&  ( 8 + 6 \beta +\frac 12 \beta^2 )  {z^2}\!\!& 
2  ( 6 + 2 \beta   )  z \!\!&\cr 
  \!\!&    \!\!&   \!\!&  \!\!&   \!\!&\!\!\!\!\ddots\cr
\end{array} \right) 
\label{ej}
\ee}}
For $\beta\equiv\bm\in{\Z}^-$ this action directly provides, by means of a
quotient space, the finite dimensional representations   of dimension
$|\bm-1|$ much in the same way as for the classical counterpart. Therefore,
we will denote
$|\bm-1|=2j_z+1\in\Z^+$ 
 being $j_z = 0,1/2,1\ldots$ Indeed it is clear that the 
expressions (\ref{eg}) contain a power series in $z$ such whose first terms
coincides with the non-deformed analogue shown in (\ref{cf}).
We  write down   the 2, 3 and 4-dimensional matrix
representations:

\noindent
{(a)} $\bm=-1$, ${\cal C}_z=3/2$, $j_z=1/2$.
\be
\jp=\left(\begin{array}{ll}
. & .   \cr 
1  & .  \cr  
\end{array}\right) 
\qquad
\jm=\left(\begin{array}{ll}
. & 1   \cr 
-\frac 14\, z^2  & z  \cr  
\end{array}\right) \qquad
\jj=\left(\begin{array}{ll}
-1 & .   \cr 
-z & 1  \cr  
\end{array}\right)  
\label{ek}
\ee

\noindent
{(b)} $\bm=-2$, ${\cal C}_z=4$, $j_z=1$.
\be
\jp=\left(\begin{array}{lll}
  . & . & . \cr 1 & . & . \cr . & {\sqrt{2}} & .
\end{array}\right) \quad
\jm=\left(\begin{array}{lll}
. & 2 & . \cr -z^2 & 2z & {\sqrt{2}} \cr
-{\sqrt{2}}\,z^3 & {\sqrt{2}}\,z^2 & 2z
\end{array}\right)  \quad
\jj=\left(\begin{array}{lll}
 -2 & . & . \cr -2z & . & . \cr -2{\sqrt{2}}\,z^2 & . & 2 
\end{array}\right)  
\label{el}
\ee

\noindent
{(c)} $\bm=-3$, ${\cal C}_z=15/2$, $j_z=3/2$.
\bea
&&\jp=\left(\begin{array}{llll}
 . & . & . & . \cr 1 & . & . & . \cr . & 
  {\sqrt{2}} & . & . \cr . & . & {\sqrt{3}} & .  
\end{array}\right) 
\quad
\jj=\left(\begin{array}{llll}
-3 & . & . & . \cr -3z & -1 & . & . \cr -2\sqrt{2}\, z^2 & -\sqrt{2}\, z& 1
   & . \cr 
-2\sqrt{6}\, z^3 &-5\sqrt{\frac 23}\, z^2 & \sqrt{3}\,z & 3  
\end{array}\right)  \cr
&&\jm=\left(\begin{array}{llll}
 . & 3 & . & . \cr -{{9 }\over 4}\, z^2 & 3\,z & 
  2\,{\sqrt{2}} & . \cr 
-{{9 }\over {2\,{\sqrt{2}}}}\,{z^3}    & 
  {{3}\over {2\,{\sqrt{2}}}}\,{z^2}& 4\,z & {\sqrt{3}} \cr 
   -3\,{\sqrt{{3\over 2}}}\,{z^4} & 
-  \frac {\sqrt{3}}{2\sqrt{2}}\,{z^3} & 
\frac {29}{4\sqrt{3}}\,{z^2}&3\, z
\end{array}\right) 
  \label{em}
\eea

By means of the general formula for the $R$-matrix (\ref{ed})  we can get
explicit solutions of the quantum Yang-Baxter equation by  substituting the
representations just found.  The computations are considerably simplified due
to the factorized form (\ref{ed}). As an example we write down  the $4\times
4$ and $9\times 9$ $R$-matrices corresponding to the above 2 and 3-dimensional
representations, respectively: 
\be
{\cal R}=\left(\begin{array}{llll}
  1 & . & . & .   \cr 
-z & 1 & . & .   \cr 
z & . & 1 & .   \cr 
z^2 & -z & z &1   \cr 
 \end{array}\right)        
\ee
\be
{\cal R}=\left(\begin{array}{lllllllll}
 1 & . & . & . & . & . & . & . & . \cr -2\,z & 1 & . & . & .
    & . & . & . & . \cr 2\,{\sqrt{2}}\,{z^2} & -2\,{\sqrt{2}}\,z & 1
    & . & . & . & . & . & . \cr 2\,z & . & . & 1 & . & . & . & . & .
   \cr . & . & . & . & 1 & . & . & . & . \cr . & 
  2\,{\sqrt{2}}\,{z^2} & -2\,z & . & . & 1 & . & . & . \cr 
  2\,{\sqrt{2}}\,{z^2} & . & . & 2\,{\sqrt{2}}\,z & . & . & 1 & . & 
  . \cr . & . & . & 2\,{\sqrt{2}}\,{z^2} & . & . & 2\,z & 1 & . \cr 
  . & -4\,{z^3} & 2\,{\sqrt{2}}\,{z^2} & 4\,{z^3} & . & 
  -2\,{\sqrt{2}}\,z & 2\,{\sqrt{2}}\,{z^2} & 2\,{\sqrt{2}}\,z & 1
  \end{array}\right)       
\ee

Upper bounded representations inside polynomial spaces
$\{x^n\}_{n=0}^{+\infty}$ are supplied by a differential (difference)
realization by means of the operator 
\be
D_z \equiv ({e^{2z\partial_x}
-1})/{2z} .
\label{eoo}
\ee
 The action of $D_z$ is just a discrete derivative:
$D_z\phi(x) = ({\phi(x+2z) - \phi(x)})/{2z}$. Thus, we have
\bea
&& \jp=\partial_x,  \qquad \jj=2D_z x + z\beta D_z +\beta,\cr
&& \jm=-D_zx^2 - z\beta  D_z x - \frac{z^2\beta^2}4 
D_z - {\beta} x.
\label{eo}
\eea
From these expressions, we see that non-standard  quantum deformations are
related to a discrete difference calculus quite different to that of standard
$U_qsl(2,\R)$. In  analogy to the classical case, the finite dimensional
representations originated from (\ref{eo}) for $\bp-2\in {\bf Z}^+$  and
$\bp-1 = 2j_z +1\in {\bf Z}^+$ is also supported by $\langle 1,x,\ldots
x^{\bp-2}\rangle$. These representations, denoted by $j_z$, will  be put to
work in some examples along the next section.

%%%%%%%%%%%COPRODUCT REPRESENTATIONS%%%%%%%%%%%%%%%%%%%%%%%%%%%%%%%%%%%%%

\sect{Tensor product representations and decomposition rules}

Given a pair of representations for the $U_zsl(2,\R)$  algebra acting on 
the vector spaces ${\cal H}_{1}$, ${\cal H}_{2}$ the  coproduct (\ref{ea}) 
originates a new representation in the tensor product  space ${\cal
H}_1\otimes {\cal H}_2$. Although the initial representations may be
irreducible the final coproduct representation  will in general be reducible.
We shall see that, for some particular finite dimensional cases worked out
below, the coproduct representation is completely reduced into irreducible
components following the same well known rules valid for the classical
$sl(2,\R)$ integer representations. Using  the conventional notation
\cite{bl}, this decomposition can be written as 
\be
j_z\otimes j_z' = |j_z+j'_z|  \oplus |j_z+j'_z-1|\oplus \ldots
\oplus |j_z-j'_z| ,
\ee
where $j_z$ and $j'_z$ are positive half-integers corresponding to the
quantum representations of dimension $2j_z+1$ or $2j'_z+1$, respectively.
However, the vector basis of the irreducible support  subspaces expressed in
terms of the original basis (i.e., the Clebsch-Gordan coefficients) become
quite different to those of the $sl(2,\R)$ Lie algebra due to extra terms
containing powers of the deformation parameter $z$.

We shall examine these features in detail for two  simple examples by making
use of the differential realizations given in (\ref{eo})  that are
particularly easy to handle for computations.
\smallskip

\noindent $\bullet$ 
$1/2 \otimes 1/2 $ representations.

The representation $j=1/2$ for the Lie algebra $sl(2,\R)$ is realized in
the polynomial vector space spanned by $\{|-1\rangle =1,|-2\rangle=x\}$.
For $j=1$ the basis is chosen in the form 
$\{|-1\rangle=1,|-2\rangle=x,|-3\rangle=x^2\}$; 
we shall use the variable
$y$ for the second space in the tensor product.   The coproduct
representation in this case is defined by
$\Delta^{\rm clas}(X) = 1\otimes U_{1/2}(X) + U_{1/2}(X) \otimes 1$, where
$X$ is any of the algebra generators and $U_{1/2}$ is for the $j=1/2$-spin
representation. The decomposition $1/2 \otimes 1/2 = 1\oplus 0$ has support
spaces whose basis are 
\bea
{\cal H}^{clas}_1 &=& \langle 
{\bf E}^{clas}_{-1}\equiv 1\, ,
{\bf E}^{clas}_{-2}\equiv \frac12(x+y) \, ,
{\bf E}^{clas}_{-3}\equiv xy \rangle,\nonumber\\
{\cal H}^{clas}_0 &=&  
\langle {\bf U}^{clas}_{-1} \equiv \frac12(x-y) \rangle .
\label{en}
\eea
Obviously, the triplet generating ${\cal H}^{clas}_1$ is  symmetric under the
permutation map $\sigma (a\otimes b)= b\otimes a$, while  the singlet
underlying ${\cal H}^{clas}_0$ is antisymmetric.

Now, for the non-standard $U_zsl(2,\R)$ the coproduct representation is
defined according to (\ref{ea}), and we consider the
value $\bp = 3$ that corresponds to the case $j_z=1/2$ (\ref{ek}). The
reduction $1/2  \otimes 1/2  = 1 \oplus 0 $ keeps on correct too. Here, the
representation
$j_z=1 $ is obtained with $\bp = 4$  
(\ref{el}), while $0$ is of course for the trivial representation. Explictly,
the invariant vector subspaces 
${\cal H}_1 = \langle {\bf E}_{-1}, {\bf E}_{-2},{\bf E}_{-3} \rangle$ and
${\cal H}_0 = \langle {\bf U}_{-1}\rangle$ are as follows  in terms of the
basis (\ref{en}):
\bea
{\bf E}_{-1} &=& {\bf E}^{clas}_{-1} \nonumber\\
{\bf E}_{-2} &=& {\bf E}^{clas}_{-2} \nonumber\\
{\bf E}_{-3} &=& {\bf E}^{clas}_{-3} + \frac{3z^2}{4} {\bf E}^{clas}_{-1} +
z{\bf U}^{clas}_{-1} \label{cob}\\
{\bf U}_{-1} &=& {\bf U}^{clas}_{-1} + \frac z2{\bf E}^{clas}_{-1} .
\nonumber
\eea
Note that symmetry in the basis of ${\cal H}_1$ and 
antisymmetry in ${\cal H}_0$ do not hold unless we  assume that  the
permutation map $\sigma$ transforms the deformation parameter $z$ into $-z$.
On the other hand, it can be easily proven that, for all $z$, the
transformation (\ref{cob}) is always well-defined. Therefore, roots of unity
seem to be not privileged for the non-standard deformation.

\noindent
$\bullet$   $1 \otimes 1/2 $ representations.

First we shall supply the basis of ${\cal H}_{1 \otimes 1/2} =
{\cal H}^{clas}_{3/2} \oplus {\cal H}^{clas}_{1/2}$ for the
reduction
$1\otimes 1/2=3/2
\oplus 1/2$ in the Lie algebra context. For the
$j=3/2$ representation we use the basis 
$\{|-1\rangle =1,|-2\rangle =x,|-3\rangle =x^2,|-4\rangle =x^3\}$,  and the
invariant subspaces in the classical tensor product are spanned by
\bea
{\cal H}^{clas}_{3/2} &=& \langle 
{\bf E}^{clas}_{-1}\equiv 1\, ,
{\bf E}^{clas}_{-2}\equiv \frac13(y +2x) \, ,
{\bf E}^{clas}_{-3}\equiv \frac13(2xy +x^2) \, ,  
{\bf E}^{clas}_{-4}\equiv x^2y
\rangle ,\cr
{\cal H}^{clas}_{1/2} &=&  
\langle 
{\bf U}^{clas}_{-1} \equiv \frac12(y-x) \, ,
{\bf U}^{clas}_{-2} \equiv \frac12{xy - x^2} .
\rangle
\eea

With respect to the deformed quantum algebra $U_zsl(2,\R)$ it can be  checked
directly that its coproduct leads to the same direct sum reduction 
$1 \otimes 1/2 =3/2  \oplus 1/2$ on the same polynomial vector space
${\cal H}_{1 \otimes 1/2}$ but with new invariant subspaces 
${\cal H}_{3/2}=
\langle {\bf E}_{-1},  {\bf E}_{-2},  {\bf E}_{-3},  {\bf E}_{-4}\rangle$ 
and 
${\cal H}_{1/2}=\langle {\bf U}_{-1}, {\bf U}_{-2}\rangle$ 
given by the following deformed change of basis:
\bea
{\bf E}_{-1} &=& {\bf E}^{clas}_{-1}\nonumber\\
{\bf E}_{-2} &=& {\bf E}^{clas}_{-2} \nonumber\\
{\bf E}_{-3} &=& {\bf E}^{clas}_{-3} +\frac{3z^2}{4}{\bf E}^{clas}_{-1} 
-\frac{2z}3{\bf U}^{clas}_{-1} \nonumber\\
{\bf E}_{-4} &=& {\bf E}^{clas}_{-4} + \frac{9z^2}{4} {\bf E}^{clas}_{-2} - 
\frac{9z^3}{4}{\bf E}^{clas}_{-1} -2z {\bf U}^{clas}_{-2} -
\frac{z^2}{3}{\bf E}^{clas}_{-1}\\
{\bf U}_{-1} &=& {\bf U}^{clas}_{-1} - \frac{z}{2}{\bf E}^{clas}_{-2} 
\nonumber\\
{\bf U}_{-2} &=& {\bf U}^{clas}_{-2}  - \frac{3z}{8}{\bf E}^{clas}_{-3} + 
\frac{3z^2}{8}{\bf E}^{clas}_{-2} . \nonumber 
\eea
As expected, the limit $z\to 0$ provides the classical partners of  the
reduction process.

At this point it is worth mentioning that representations  of the standard
deformation of $sl(2,\R)$ are strongly different from their non-standard
counterparts. On one hand, such standard representations can be essentially
constructed by substituting some matrix elements of the classical matrices by
the corresponding $q$-numbers \cite{bibook} and, consequently, the same holds
for the Clebsch-Gordan coefficients. This straightforward method is no longer
valid for the non-standard case, where $q$-numbers do not work and, moreover,
some new non-vanishing Clebsch-Gordan coefficients have to be added  with
respect to the classical theory. Since $q$-numbers  are directly related to
the peculiarities af roots of unity, the lost of such properties in  the
non-standard case seems quite natural.

%%%%%%%%%%%%%%%%%%%%%%%%%%%%%%%%%%%%%%%%%%%%%%%%%%
%%%%%%%%%%%%%%%%%Representaciones a dos bosones%%%%%%%%%%%

\sect{Two-boson $U_zsl(2,\R)$ representations and their contraction to
Poincar\'e}

As it has been shown above  the
one-boson representations of   $sl(2,\R)$ are
closely linked with its interpretation as a one-dimensional conformal
algebra. In contrast,   a description in
terms of   two-boson algebras is physically related to its
role as  a (1+1)-dimensional kinematical algebra.  This fact allows us to
perform a contraction in order to reach the  (1+1) Poincar\'e algebra
representations; such a process  cannot be applied  onto the one-boson
representations of section 3.

\subsect{Classical two-boson  representations}

Let us again begin with a discussion for $sl(2,\R)$ since it will give
us a natural reference. We consider two independent boson
algebras 
\be
[\aam,\aap]=1,\qquad [\bbm,\bbp]=1.
\label{fa}
\ee
   A two-boson representation of $sl(2,\R)$ is the
following
\bea
&&\jp=\aap,\qquad 
 \jj=2\aap\aam - 2\bbp\bbm,\cr
&&\jm=-\aap\aam^2 + 2\bbp\bbm\aam +\alpha\bbp.
\label{fb}
\eea
There are, of
course a high arbitrariness for many other expressions realizing the
same $sl(2,\R)$ algebra, but the simple choice shown in (\ref{fb}) will be
enough for our purposes. The only track left here
of such wide range of possibilities is the free parameter
$\alpha$. This parameter should not be seen as an irreducible representation
label as it was the case of $\beta$, for instance, with respect to
(\ref{cd}). We will call formulas (\ref{fb}) the two-boson GD
realization. They should be compared to the more common Jordan--Schwinger
realization  \cite{sc}.

A differential realization is obtained by taking 
\be 
 \aam=-x^+  ,\qquad \aap=\frac{\partial}{\partial x^+}\equiv
\partial_+,\qquad
 \bbm=-x^- ,\qquad \bbp=\frac{\partial}{\partial x^-}\equiv \partial_-, 
\label{fc}
\ee 
where, $x^+=t+x$ and $x^-=t-x$ can be identified   as light-cone
coordinates. Then, we get
\bea
&&\jp=\partial_+,\qquad 
 \jj=-2x^+\partial_+  +2x^-\partial_-,\cr
&&\jm=-(x^+)^2\partial_+ + 2 x^+x^-\partial_-  +\alpha \partial_- .
\label{fd}
\eea
The expression for the Casimir (\ref{cb})  acting on the wavefunction space
$\phi(x^+,x^-)$ is the second order operator 
\be
{\cal C} = 2(x^-)^2\partial_{- -}^2 + 4 x^-\partial_- + 2\alpha 
\partial_{+-}^2 .
\ee 

 By means of the following  In\"on\"u--Wigner contraction, to be applied  in
(\ref{ca}),  
\be
P_+=\varepsilon J_+,\qquad P_-=\varepsilon J_-,\qquad K=\frac 12 J_3,
\label{fh}
\ee
the boost $K$ and the light-cone traslations $P_\pm$ of the (1+1)  Poincar\'e
algebra ${\cal P}(1+1)$ are generated in the limit $\varepsilon\to 0$, since
the comutation rules between these new generators are
\be
[K,P_+]=P_+,\qquad [K,P_-]=-P_-,\qquad [P_+,P_-]=0.
\label{fi}
\ee

In order to contract the representation  (\ref{fb}) we  consider (\ref{fh})
together with 
\cite{Marc}: 
\be  
\aam\to \varepsilon^{-1}\aam,\qquad \aap\to
\varepsilon\aap,\qquad 
\bbm\to  \bbm, \qquad \bbp\to  \bbp,\qquad
\alpha\to \varepsilon\alpha  ,
\label{fk}
\ee 
and the limit $\varepsilon\to 0$ provides the two-boson ${\cal P}(1+1)$
representations
\be 
 P_+=\aap,\qquad
 K= \aap\aam  -  \bbp\bbm,\qquad
 P_-= \alpha \bbp .
\label{fl}
\ee

\subsect{Quantum two-boson  representations: the Poincar\'e algebra}

As far as the quantum non-standard algebra $U_zsl(2,\R)$ is concerned, the
corresponding two-boson version takes the form
\bea
&& \jp=\aap,\qquad 
  \jj=\frac{e^{2 z \aap} -1  }{z}\,\aam 
- 2\bbp\bbm,\cr
&& \jm=-\frac{e^{2 z \aap} -1  }{2z}\,\aam^2 + 2\bbp\bbm\aam + \alpha 
\bbp+2z( \bbp\bbm +  \bbp^2\bbm^2) .
\label{fg}
\eea

Algebraically $U_zsl(2,\R)$ can be contracted to a non-standard
$U_z{\cal P}(1+1)$ algebra by defining the generators as in
(\ref{fh}) and at the same time setting 
\be
z\to \varepsilon^{-1}z ,
\label{ffoo}
\ee
 so that we get the Hopf algebra
\bea
&&  \Delta (P_+) =1 \otimes P_+  + P_+ \otimes 1,\cr
&& \Delta (K) =1 \otimes K  + K\otimes
e^{2z P_+ }, \label{fn} \\
&& \Delta (P_-) = 1 \otimes P_- + P_-\otimes
e^{2 z P_+ }, 
\nonumber
\eea
\be
\epsilon(X) =0,\qquad 
 \mbox{for $X\in \{K,P_+,P_-\}$},
 \ee
\be
\gamma(P_+)=-P_+,\qquad \gamma(K)=-K e^{-2zP_+},\qquad 
\gamma(P_-)=-P_- e^{-2zP_+},
\ee
\be 
[K,P_+ ]= \frac{e^{2 z P_+} -1  } {2z},\qquad 
[K,P_-]=-  P_-  ,\qquad  [P_+  ,P_- ]=0.
\label{fo} 
\ee 
The quantum Casimir is found   by
contracting (\ref{ec}) as the $\lim_{\varepsilon\to 0}(\varepsilon^2 {\cal
C}_z)$  and the $R$-matrix comes directly from the contraction of
(\ref{ed}):
\bea
&&{C}_z=  \frac {1-e^{-2z P_+}}{ z}\,P_- \\
&&{R}=\exp\{- 2z P_+\otimes K\}\exp\{2z K\otimes P_+\} .
\eea
The corresponding classical $r$-matrix is $r= 2 z K\wedge P_+$.
These results are in full concordance with those obtained
in \cite{tt} by a $T$-matrix approach.

After applying  the transformations (\ref{fh}), (\ref{fk}) and (\ref{ffoo})
on (\ref{fg}), the contracted two-boson representation becomes
\be 
  P_+=\aap,\qquad
  K=\frac{e^{2 z \aap} -1  }{2z}\,\aam 
- \bbp\bbm,\qquad
  P_-=  \alpha \bbp
 .\label{fp}
\ee

%\newpage

%%%%%%%%%%%%%%%%%%%%%%%%%%%%%%%%%%%%%%%%%%%%%%%%
%%%%%%%%%%%%%%%%%La extension del no estandar%%%%%%%%%%%

\sect{Extended $U_zsl(2,\R)$ and its contraction to the 
oscillator algebra}

The infinite dimensional representations for the non-standard quantum
oscillator algebra \cite{osc} can be deduced by performing a contraction 
on a (non-standard) quantum deformation $U_zsl(2,\R)$ of  the pseudo-extended
$sl(2,\R)$ Lie algebra. In this section we develop such a process at both
classical and quantum levels.

\subsect{Classical level}

It is well-known that a trivial central extension of the Lie algebra
$sl(2,\R)$ leads through a careful contraction to a non-trivial
extension of  the (1+1) Poincar\'e algebra corresponding to a constant
non-null  background field \cite{no}; this contracted extended algebra is
isomorphic to the oscillator $h_4$ Lie algebra. Let us review here such
properties for the sake of comprehension and unification of notation. The
trivial extension, designed by
$\overline{sl}(2,\R)$, obeys to the commutation rules
\be
[\jj,\jp]=2\jp,\qquad  [\jj,\jm]=-2\jm,\qquad [\jp,\jm]=\jj - I ,
\qquad [I,\cdot\,]=0,
\label{ga}
\ee
where $I$ is the central extension generator. The second order Casimir is now
\be
{\cal C}=\frac 12 \jj^2 - \jj I +\jp\jm+\jm\jp .
\label{gb}
\ee

The one-boson realization for (\ref{ga}) is
\be
\jp = a_+,\quad
\jj = 2a_+a_- + \beta 1,\quad
\jm = - a_+a_-^2 - \beta \aam + \delta \aam ,\quad
I=\delta 1,
\label{gc}
\ee
where $\delta$ and $\beta$ are free parameters related with the eigenvalue
of the Casimir by
\be
{\cal C}=\beta(\beta/2 - 1) + \delta(1-\beta) .
\label{gd}
\ee

An In\"on\"u-Wigner
contraction can be applied by defining the new generators
\be
\ap = \varepsilon \jp,\quad
\am = \varepsilon \jm,\quad
N = \jj/2,\quad
M = \varepsilon^2 I,
\label{ge}
\ee
so that  in the limit $\varepsilon \to 0$ we reach the oscillator
$h_4$ Lie algebra,
\be
[N,A_+] =  A_+,\qquad [N,A_-] = -  A_-,\qquad [A_-,A_+]=M,\qquad
[M,\cdot\,] =0.
\label{gf}
\ee
The corresponding second order Casimir is obtained as
$\lim_{\varepsilon \to 0}( - \varepsilon^2 {\cal C})$:
\be
{C}=2N M - A_+ A_- - A_- A_+.
\label{gg}
\ee 

The additional replacements 
\be
 \aam
\to \varepsilon^{-1}\aam,\qquad \aap \to \varepsilon \aap ,
\qquad \beta\to  \beta/2,
\qquad \delta \to \varepsilon^2 \delta ,
\label{gh}
\ee
provide  the  one-boson  $h_4$ realization: 
\be
N= \aap\aam + \beta,\qquad \ap = \aap,\qquad \am=\delta\aam,
\qquad M=\delta 1.
\label{gi}
\ee
Hence the eigenvalue of (\ref{gg}) is $C= \delta(2\beta-1)$.

As for the infinite dimensional irreducible representations of $h_4$ in
number state spaces, they are trivially derived from those of the boson
algebra (\ref{cdd}) and (\ref{cgg}). Note that (\ref{gi})  clarifies the
difference between considering a boson algebra and the harmonic oscillator
algebra $h_4$.

\subsect{Quantum level: the non-standard oscillator}

In the quantum context following closely the classical
approach first we must define an appropriate  extension of $U_zsl(2,\R)$ by
the addition of a new central generator $I$. In this way we will call the
extended  non-standard quantum algebra of ${sl}(2,\R)$  
 to the Hopf algebra denoted $U_z\overline{sl}(2,\R)$ and given by
\bea
&&  \Delta (\jp) =1 \otimes \jp  + \jp \otimes 1,\cr
&& \Delta (\jj) =1 \otimes \jj  + \jj\otimes
e^{2z \jp }, \label{gj} \\
&& \Delta (\jm) = 1 \otimes \jm + \jm\otimes e^{2 z \jp } + z \jj\otimes I\,
e^{2 z \jp }, \cr
&&  \Delta (I) =1 \otimes I  + I \otimes 1,\nonumber
\eea
\be
\epsilon(X) =0,\qquad 
 \mbox{for $X\in \{\jj,\jp,\jm,I\}$},
 \ee
\bea
&&\gamma(\jp)=-\jp,\qquad \gamma(I)=-I,\qquad  \gamma(\jj)=-\jj
e^{-2z\jp}, \cr
&&\gamma(\jm)=-\jm e^{-2z\jp} + z \jj I\, e^{-2z\jp},
\eea
\bea 
&&[\jj,\jp ]= \frac{e^{2 z \jp} -1  } z,\qquad 
[\jj,\jm]=-2 \jm +z \jj^2,\cr
&&[\jp  ,\jm ]= \jj - I\,e^{2z\jp},\qquad [I,\cdot\,]=0.
\label{gk} 
\eea 

The (coboundary) Lie bialgebra underlying this Hopf algebra is again
generated by  $r= z\jj\wedge\jp$.
Note that the new generator $I$ remains central and primitive;   there is
another  quantum Casimir given by
\be
{\cal  C}_z=\frac 12 \jj\, e^{-2z\jp}\jj - \jj I + \frac
{1-e^{-2z\jp}}{2z}\,\jm +\jm\,\frac {1-e^{-2z\jp}}{2z} + e^{-2z\jp}-1.
\label{ggk}
\ee

The Hopf subalgebra generated by $\jj$ and $\jp$ is the same as in the
non-extended case, therefore, the universal $R$-matrix (\ref{ed}) is 
obviously a solution of the quantum YBE for $U_z\overline{sl}(2,\R)$.
Furthermore, cumbersome computations show this $R$-matrix   
verifies ${\cal R}\Delta(\jm){\cal R}^{-1}=\sigma\circ \Delta(\jm)$ in
$U_z\overline{sl}(2,\R)$ (the proof for $I$ is trivial),  hence we conclude
that element is also a universal
$R$-matrix for the whole $U_z\overline{sl}(2,\R)$. 
At this point it is worth
mentioning that another quantum deformation of the extended 
$\overline{sl}(2,\R)$ has been recently proposed in \cite{abb} leading to a
deformed oscillator algebra with classical commutation relations, but not
keeping the aforementioned subalgebra.

The one-boson realization of $U_z\overline{sl}(2,\R)$ turns out to be a
slight modification  of the non-extended case (\ref{ee}). Besides the new
generator $I = \delta 1$, the only generator that changes is $J_-$ that reads
\be
\jm=-\frac{e^{2 z \aap} -1  }{2\,z}\,\aam^2 - 
\beta \, \frac{e^{2 z \aap} +1  }{2}\,\aam 
- z \beta^2\,\frac{e^{2 z \aap} -1  }{8} + \delta e^{2z\aap}\aam + 
\frac{\beta z}{2} \delta e^{2z\aap} .
\label{gl}
\ee

Now we proceed to carry out the contraction from $U_z\overline{sl}(2,\R)$
to the non-standard quantum oscillator algebra \cite{osc}, denoted   $U_z
h_4$. At the Hopf algebra level, we consider the new generators defined by
(\ref{ge}) and also the transformation of the deformation parameter 
$z$  (\ref{ffoo}).   Thus, when $\varepsilon \to 0$ we arrive at the Hopf
structure of $U_zh_4$ given by
\bea
&&\Delta(\ap)=1\otimes \ap +\ap \otimes 1,\qquad 
 \Delta(\bb)=1\otimes \bb +\bb \otimes 1,\cr
&&\Delta(\aa)=1\otimes \aa +\aa \otimes e^{2z\ap},\cr 
&& \Delta(\am)=1\otimes \am +\am \otimes e^{2z\ap}+2z\aa\otimes \bb\,
e^{2z\ap}, \label{gn}
\eea
\be
\epsilon(X)=0,\qquad X\in\{\aa,\ap,\am,\bb\},
\label{go}
\ee
\bea
&&\gamma(\ap)=-\ap,\qquad \gamma(\bb)=-\bb,\cr 
&&\gamma(\aa)=-\aa  e^{-2z\ap},\quad 
\gamma(\am)=-\am  e^{-2z\ap}+ 2 z \aa\bb  e^{-2z\ap}, 
\label{gp}
\eea
satisfying the commutators
\be 
 [\aa,\ap]=\frac{e^{2z\ap}-1}{2z},\quad [\aa,\am]=-\am
,\quad [\am,\ap]=\bb  e^{2z\ap},\quad 
[\bb,\cdot\,]=0 ,
\label{gq}
\ee 
where the classical $r$-matrix is $r=2z \aa\wedge \ap$.
Besides the generator $\bb$  there is another central operator
 which is directly obtained from (\ref{ggk}) by means of the
$\lim_{\varepsilon\to 0} (-\varepsilon^2 {\cal C}_z)$ giving
rise to the expression
\be
{C}_z=2\aa\bb+  
\frac{e^{-2z\ap}-1}{2z}\, \am +\am \, \frac{e^{-2z\ap}-1}{2z}.
\label{gr}
\ee
Likewise, the corresponding  universal $R$-matrix is found by contracting
(\ref{ed}):
\be
R=\exp\{-2z\ap\otimes\aa\}\exp\{2z\aa\otimes\ap\} .
\label{gs}
\ee

The boson representation of $U_z(h_4)$  can be found using the same routine
taking into account (\ref{ge}), (\ref{ffoo}) plus the extra replacements
(\ref{gh}).  Thus we have
\bea
&& \ap=\aap,\qquad  \bb=\delta 1 ,\cr
&& \am=\delta e^{2z\aap}\aam +  \delta{\beta z}\,  e^{2z\aap} 
,\label{gtq}\\ 
&& \aa= \frac{e^{2z\aap}-1}{2z}\,\aam+\beta \,
\frac{e^{2z\aap}+1}{2} .\nonumber
\eea
Hence the action  on the states
$\{|m\rangle\}_{m=0}^{\infty}$  reads
\bea
&&\ap|m\rangle =\sqrt{m+1}\, |m+1\rangle ,\qquad
M|m\rangle =\delta\, |m\rangle , \cr
&&\am|m\rangle =\delta\,\sqrt{m}\, |m-1\rangle + \delta  \sum_{k=0}^\infty
\frac{(2z)^{k+1} }{k!} \, \sqrt{\frac {(m+k)!}{m!} }\left(
\frac{m}{k+1}+\frac{\beta}{2}\right)\, |m+k\rangle ,\cr
&&\aa|m\rangle =(m+\beta)\, |m\rangle +  \sum_{k=1}^\infty \frac{(2z)^k
}{k!} \, \sqrt{\frac {(m+k)!}{m!} }\left(
\frac{m}{k+1}+\frac{\beta}{2}\right)\, |m+k\rangle  .\label{gt} 
\eea
The explicit   matrix form is
\be
 \am=\delta\left (\begin{array}{llllll}
\beta z  & 1 & . & . & .& \cr 
2 \beta  z^2 &    (2 + \beta)z  & \sqrt{2} & . & .& \cr 
2{\sqrt{2}} \beta  {z^3} &
  2{\sqrt{2}}  ( 1 + \beta   )  z^2 & (4 +
\beta) z   & \sqrt{3} & . &\cr 
\frac {8}{\sqrt{6}} \beta  {z^4} & 
   2 {\sqrt{6}}  ( {2\over 3} + 
 \beta  )  {z^3} & 
 2 {\sqrt{3}}  ( 2 + \beta   )  z^2 & (6 + \beta) z 
    & \sqrt{4} &\cr 
\frac {8}{\sqrt{6}} \beta  {z^5} & 
  8 {\sqrt{{2\over 3}}}  ( \frac 12 + \beta   )  
   {z^4} & 4 {\sqrt{3}} 
    ( {4\over 3} + \beta   )  {z^3} & 
   4   ( 3+ \beta   )  z^2 & (8 + \beta) z & \cr
  &    &   &  &   &\ddots\cr
\end{array} \right),
\nonumber
\ee
\be
\aa=\left (\begin{array}{llllll}
\beta  & . & . & . & .& \cr 
\beta  z &    1 + \beta  & . & . & .& \cr 
{\sqrt{2}} \beta  {z^2} &
  {\sqrt{2}}  ( 1 + \beta   )  z & 2 +
\beta   & . & . &\cr 
\frac {4}{\sqrt{6}} \beta  {z^3} & 
   2 {\sqrt{6}}  ( {1\over 3} + 
  \frac 12   \beta  )  {z^2} & 
  {\sqrt{3}}  ( 2 + \beta   )  z & 3 + \beta 
    & . &\cr 
\frac {4}{\sqrt{6}} \beta  {z^4} & 
  4 {\sqrt{{2\over 3}}}  ( \frac 12 + \beta   )  
   {z^3} & 2 {\sqrt{3}} 
    ( {4\over 3} + \beta   )  {z^2} & 
   2  ( 3 + \beta   )  z & 4 + \beta & \cr
  &    &   &  &   &\!\!\!\!\ddots\cr
\end{array} \right), 
  \label{gu}
\ee 
being (\ref{ci}) the matrix for $\ap$ and $\delta {\bf 1}$ with ${\bf 1}$ the
identity matrix for $M$.

Finally,  the explicit infinite dimensional
representations in the monomial basis
$\{x^n\}$  as well as a differential difference 
realization in terms of the operator (\ref{eoo}) can be readily obtained.
The latter is: 
\bea
&& \ap=\partial_x, \qquad M=\delta 1, \qquad \aa= D_z x + z\beta D_z
+\beta,\cr 
&& \am=2\delta z D_zx  + 2 z\delta \beta z^2 D_z  + \delta x +
\delta\beta z.
\label{gv}
\eea

\sect{Quadratic algebras}

It is sometimes convenient to choose a suitable basis for the quantum Hopf
algebras in order to get quadratic  commutators.  In doing so the computation
of commutators and representations at the algebra level is greatly
facilitated, although this procedure  can add some   extra complications at
the coalgebra level. We shall briefly consider here that point with the help
of an auxiliary deformed boson algebra:
\be
 \baap =\frac{e^{2z\aap}-1}{2z},\qquad \baam= \aam +\mu z ,
\label{sa}
\ee
satisfying
\be
[\baam,\baap] = 1+2z\baap,
\label{ha}
\ee
where $\mu$ is a free parameter to be fixed for each case.

\noindent
{\it (a) Quadratic $U_zsl(2,\R)$}. 

Let us change to the basis
\be
\bjp = \frac{e^{2z\jp}-1}{2z},\quad 
\bjj = \jj,\quad 
\bjm = \jm
\label{hb}.
\ee
Then, a (deformed) boson realization is written  in the following form by
taking $\mu= \beta/2$: 
\be
\bjp = \baap,\quad 
\bjj = 2\baap\baam +\beta 1,\quad 
\bjm = -\baap\baam^2 - \beta\baam + \frac12 \beta^2z,
\label{hc}
\ee
to be confronted with (\ref{cd}). The commutators are easily derived from
(\ref{ha}):
\bea
&&[\bjj,\bjp] = 2\bjp(1+2z\bjp),\nonumber\\
&&[\bjj,\bjm] = -2\bjm +z\bjj^2,\\
&&[\bjp,\bjm] = (\bjj -2z\bjp)(1+2z\bjp).\nonumber
\eea

\noindent
{\it (b) Quadratic non-standard Poincar\'e}.

Here we can choose $\mu =0$, i.e., $\baam=\aam$. Let us define the basis
\be
\overline{P}_+ = \frac{e^{2zP_+}-1}{2z},\quad
\overline{P}_-  = P_-,\quad
\overline{K} =  K.
\ee
The boson realization is given in terms of  $\{\baap,\baam\}$ and $\{
b_+,b_-\}$,
\be
\overline{P}_+ = \baap,\quad
\overline{P}_-  =\alpha b_+,\quad
\overline{K} =  \baap\baam - b_+b_-,
\ee
cf. equations (\ref{fl}). Now the commutators are
\be 
 [\overline{K},\overline{P}_+] = \overline{P}_+(1 +2z\overline{P}_+) 
,\quad
 [\overline{K},\overline{P}_-] =  \overline{P}_-,\quad
 [\overline{P}_+,\overline{P}_-] = 0. 
\ee

\noindent
{\it (c) Quadratic $U_zh_4$}. 

Now we take the basis generators in the form
\be
\bap = \frac{e^{2z\ap}-1}{2z},\quad
\bam = e^{-2z\ap}\am,\quad
\baa = \aa,\quad
\overline{M} = M .
\ee
The corresponding (deformed) boson realization for $\mu=\beta$ reads
\be
\bap = \baap,\quad
\bam = \delta\baam,\quad
\baa =  \baap\baam + \beta 1,\quad
\overline{M} = \delta 1   .
\ee
to be compared with (\ref{gi}). The commution rules become
\bea
&&[\baa,\bap] =  (1+2z\bap)\bap,\quad [\baa,\bam] = -(1+2z\bap)\bam,\cr
&& [\bam,\bap] = (1+2z\bap)\overline{M},\quad [\overline{M},\cdot\,]=0 .  
\eea

%%%%%%%%%%%%%%%%% Concluding remarks %%%%%%%%%%%

\sect {Concluding remarks}

Along this paper we have given a unified treatment for a class of
non-standard algebras related to $U_zsl(2,\R)$: the extended
$U_z\overline{sl}(2,\R)$, Poincar\'e $U_z{\cal P}(1+1)$, and oscillator
$U_zh_4$ deformed algebras. All these algebras share the same Hopf subalgebra
(in the $U_zsl(2,\R)$ case is generated by $\{\jj,\jp\}$) which  leads to a
common universal $R$-matrix.

At the same time we have computed the representations of these non-standard
algebras by means of boson operators. We have remarked when it was suitable
the use of one or two-boson algebras to describe each type of
representations. The general conclusion is that the functions involved
generalize the GD map for the angular momentum, in contrast to
the usual Jordan--Schwinger map, that turns out to be more adequate  for
standard deformations. We have obtained simple closed  expressions and
practical differential realizations for the $U_zsl(2,\R)$ representations
(see also in this respect
\cite{ab}) which parallel the Lie algebraic classification. For instance, 
finite dimensional representations have been shown to be labelled by integers
$j_z$ and we have proven how their coproduct representations  decompose by
following exactly the classical addition of angular momenta \cite{bl}.

Contraction processes onto the (one or two) bosonic
representations relating all these non-isomorphic quantum  algebras have been
found. A careful attention was paid to define the most adequate non-standard
quantum deformation of the centrally pseudo-extended algebra
$\overline{sl}(2,\R)$. Indeed, it has some original features  with respect to 
other extensions already defined in the context of  quantum deformations; for
example the altered commutator was not by writing a combination of the
primitive generators as usual. At the same time the coproduct of the initial
generators was also changed with the help of the new primitive generator $M$.
This extension allows for a contraction to $U_zh_4$ keeping intact the Hopf
subalgebra $\{\jj,\jp\}$.

We would also like  to stress that we have preferred to present the
properties under study in parallel with the classical situation at the Lie
algebra level. In this way we wanted to keep closer to the physical language
as well as to appreciate the prime role played by
GD like maps in non-standard deformations.

%%%%%%%%%%%%%%%%% ACKNOWLEDGEMENTS %%%%%%%%%%%

%\newpage

\bigskip
\bigskip

\noindent
{\large{{\bf Acknowledgements}}}

\bigskip

 This work has been
partially supported by DGICYT (Projects PB94--1115 and PB95--0719) from the
Ministerio de Educaci\'on y Ciencia de Espa\~na.

\bigskip
\bigskip

%%%%%%%%%%%%%%%%% BIBLIOGRAPHY %%%%%%%%%%%
%\footnotesize

\end{document}